\documentclass[11pt]{amsart}
\usepackage{fullpage, amsfonts,amsmath,amsthm,amssymb,hyperref,float,mathtools,color, latexsym}

\title{Symplectic invariants of semitoric systems and the inverse problem for quantum systems}
\date{}
\author{{\'A}lvaro Pelayo}

\begin{document}
\maketitle

\begin{center}
\emph{In memory of  Professor Johannes (Hans) J. Duistermaat (1942-2010)}
\end{center}

\begin{abstract}
 Simple semitoric systems were classified about ten years ago 
  in terms of a collection of invariants, essentially given by a convex polygon with 
 some marked points corresponding to focus\--focus singularities. 
 Each marked point is endowed with labels which are symplectic
  invariants of the system.  We will review the construction of these
   invariants, and explain how they have been  generalized 
  or applied in different contexts. One of these applications concerns 
  quantum integrable
systems and the corresponding inverse problem, which 
asks how much information of the associated classical 
system can be found in
the spectrum. An approach to this problem has been 
 to try to compute invariants in the spectrum. We will explain how 
 this has been recently achieved  for some of the  invariants of semitoric 
 systems, and discuss an open question in this direction. 
  \end{abstract}

\section{Introduction}

An integrable system
on a symplectic manifold $(M,\omega)$ of dimension $2n$ is given by $n$ smooth functions 
$$f_1,\ldots, f_n \colon (M,\omega) \to \mathbb{R}$$
which are independent and in involution, with respect to the symplectic form $\omega$. For example,
the spherical pendulum, the coupled angular momenta, and the Jaynes\--Cummings coupled spin\--oscillator model, are all integrable systems. 

Much effort has been made in recent years to develop a 
general symplectic theory for integrable systems, by  constructing as many 
symplectic invariants of these systems as possible.  
In a few cases a complete set of invariants has been found, which classifies all integrable systems
in a certain class. Two such classes are the class of integrable systems of toric type (on compact 
manifolds of any dimension), and the class of integrable systems of semitoric type (on $4$\--manifolds, compact or not), under some conditions. 
The first class concerns systems for which all the functions $f_1,\ldots,f_n$ generate
periodic flows of the same period, while the second class concerns systems $f_1,f_2$ for which $f_1$
generates a periodic flow but there is no requirement on $f_2$. 

 In this paper we will highlight some classical results about toric integrable systems and a few recent results about semitoric  integrable systems.  In order to do this we first introduce
the  language of symplectic geometry and integrable systems, and the main concepts one needs
to  understand the aforementioned results concerning toric and semitoric integrable systems.  Then we will
discuss quantum integrable systems. These are given by a collection of semiclassical
operators $$P^1=(P^1_{\hbar})_{\hbar \in I}, \ldots, P^n=(P^n_{\hbar})_{\hbar \in I}$$ 
on a sequence of Hilbert spaces $(\mathcal{H}_\hbar)_{\hbar \in I}$, whose principal symbols form a classical integrable system in the sense
described above. Here $I$ is a subset of $(0,1]$ which accumulates at $0$. 

We will  concentrate on quantum integrable systems  of toric and semitoric type, that is, those
whose principal symbols form an integrable system of toric and semitoric type, respectively. We
will explain why understanding the symplectic geometry of their classical counterparts plays a key role 
in the study of inverse problems for quantum integrable systems. 

 In this direction, we will discuss 
 progress on the inverse spectral conjecture for semitoric
systems~\cite[Conjecture 9.1]{PeVN11} from about ten years ago.  The conjecture essentially says that one can construct, from the data given by  the semiclassical joint spectrum of a quantum semitoric integrable system  $P_1,P_2$ (given for each fixed value of the $\hbar$ by the support of
 the joint spectral measure), the associated classical integrable system given by the principal symbols
 $f_1,f_2$ of the semiclassical operators $P_1,P_2$.    This can be illustrated with a diagram,
   where the answer to the question on the second arrow is \emph{yes} for quantum (simple)  
    semitoric systems according to the conjecture:
    
     $$
 \textup{{Quantum system}} \rightsquigarrow \textup{{Semiclassical Spectrum}} \underbrace{\rightsquigarrow}_{\textup{possible?}}  \textup{Classical system}
$$

The question of how much information is lost in 
the step of going from a quantum integrable system to its semiclassical spectrum has received 
significant attention in recent years.  
The answer to the question in the second arrow of the diagram has been shown to be 
\emph{yes} in some  cases, for example for quantum toric systems~\cite{CPVN13}. 
In the case of quantum semitoric systems, under some conditions, one can 
read off from the spectrum some properties of
the classical system~\cite{LFPeVN16}. In order to do this one needs to
understand  the symplectic invariants
of semitoric systems, and then how to compute them
in the spectrum (using, for example, microlocal analysis). We will discuss these results and
an open question in the last section of the paper,
while the previous sections are devoted to reviewing some the main concepts of 
symplectic geometry of  integrable systems which we need for the last section. The open
question concerns spectral theory of non\--simple semitoric systems; their classical counterparts
have been recently classified~\cite{PPT19}.  Throughout the paper we  often do not present the most general definitions or results, and try to convey  the main ideas instead.  

\section{Symplectic geometry of classical integrable systems}

 The term \emph{symplectic} was introduced in Weyl's book~\cite{Weyl} as an 
 analogue to \emph{complex}.  Symplectic geometry has
its roots in classical mechanics.  One can
trace the first steps to the seventeenth century, to the works of Galileo Galilei, Christiaan
Huygens, and Isaac Newton.  The phase space of a mechanical system is  modeled by a 
symplectic manifold.  
The first example of a symplectic manifold was given in 1808 by Joseph-Louis 
Lagrange in his  works on motions of planets~\cite{Lagrange1,Lagrange2}.    

After Lagrange,
two precursors of symplectic geometry were the works of Carl Gustav Jacob Jacobi and William Rowan Hamilton. Hamilton gave a deep formulation of Lagrangian mechanics around 1835. Their  
methods and ideas were influential in the modern view point, which starts with important contributions by many authors around the early 1970s, see~\cite{BAMS17,Weinstein} and the references therein.   After these developments several aspects of the subject became subjects on their own right. One these aspects concerned symplectic geometry of finite dimensional integrable systems. Next we introduce these systems and discuss
some of their properties.

\subsection{Phase space and symplectic forms}

A \emph{symplectic manifold} is a pair $(M,\omega)$ where $M$ is a smooth manifold and
$\omega$ is a non\--degenerate closed $2$\--form $\omega \in \Omega^2(M)$, called a 
\emph{symplectic form}. 
The non\--degeneracy condition on $\omega$ means that, at  every  $m \in M$, the skew\--symmetric bilinear form $\omega_m \colon {\rm T}_mM \times {\rm T}_mM \to \mathbb{R}$ is non\--degenerate.
This is a linear algebra condition which implies that $M$ is even\--dimensional
and also orientable, since $\omega^n$ defines a volume form on $M$.  The closedness of
$\omega$ implies that if $M$ is compact the cohomology class $[\omega^k]$, $k \in \{1,\ldots, n\}$, is non trivial,
so the even dimensional de Rham cohomology groups ${\rm H}_{\rm dR}^{2n}(M)$ of $M$, $k \in
\{1,\ldots, n\}$, are non trivial. 

Putting these restrictions together we conclude that the only sphere $S^n$ which
admits a symplectic form (in fact, many) is $S^2$. Other typical examples of symplectic manifolds are
  $(\mathbb{R}^{2n} ,\sum_{i=1}^n{\rm d}x_i \wedge {\rm d}y_i)$, where
  $(x_1,\ldots, x_n, y_1,\ldots, y_n)$ are coordinates on $\mathbb{R}^{2n}$,  
  the cotangent bundles $({\rm T}^*X, \sum_{i=1}^n{\rm d}x_i \wedge {\rm d}\xi_i)$
where $(x_1,\ldots,x_n)$ are coordinates on an $n$\--dimensional compact manifold $X$ and 
$(\xi_1,\ldots, \xi_n)$ are the usual cotangent conjugate coordinates, or 
$(S,\omega)$ where $S$ is a surface and $\omega$ an area form on it.

 In fact, 
$\mathbb{R}^{2n}$ is the local model for all symplectic manifolds of dimension $2n$: Darboux proved~\cite{darboux} that they are all  locally diffeomorphic to 
$(\mathbb{R}^{2n}, \sum_{i=1}^n{\rm d}x_i \wedge {\rm d}y_i)$,
so other than its dimension, symplectic manifolds have no local invariants.

Unless otherwise stated, in this paper \emph{all symplectic manifolds are assumed
to be connected}.

\subsection{Integrability of Hamiltonians}
  
Given a smooth function $f \colon (M,\omega) \to \mathbb{R}$ on a symplectic manifold there exists
a unique smooth vector field $\mathcal{X}_f$ such that
$\omega\,(\mathcal{X}_{f},\cdot)=-{\rm d}f;$ the vector field $\mathcal{X}_f$ is usually called the \emph{Hamiltonian
vector field induced by $f$}. Also often one refers to $f$ as a \emph{Hamiltonian function},
or simply a \emph{Hamiltonian}. 

If $(M,\omega)$ has dimension $2n$, a \emph{classical integrable system
 on $M$} is given by a collection of $n$ real\--valued smooth functions $f_1,\ldots, f_n$ on $M$ 
 such that: (1) $f_1,\ldots,f_n$  are in involution, that is, $\{f_i,f_j\}:=\omega(\mathcal{X}_{f_i},\mathcal{X}_{f_j})=0$ for all $i,j$ (another  way to say this is that $f_i$ is constant along the flow of  $\mathcal{X}_{f_j}$ for all $i,j$) and, (2)  $f_1,\ldots,f_n$      are independent, that is,   $\mathcal{X}_{f_1},\ldots, \mathcal{X}_{f_n}$ are linearly independent almost  everywhere in $M$.
  
 The term ``integrable system" comes from
 considering only the Hamiltonian $f_1$, and then looking for the maximal possible number of ``integrals
 of $f_1$"  (i.e. functions $f_j$ with $\{f_1,f_j\}=0$) which are independent. If $M$ has dimension $2n$, 
 one can have at most $n-1$ such integrals $f_2,\ldots, f_n$, and then the modern
 view point has been to call the joint map $F=(f_1,\ldots,f_n) \colon M \to \mathbb{R}^n$ 
 the \emph{momentum map} of the integrable system, or even more often the \emph{integrable system} itself. This terminology is 
in part motivated by Hamiltonian group actions where  $F$ is
the momentum of the action~\cite{[62],[95]}.

\subsection{Example of an integrable system: the spherical pendulum} \label{nm}

The spherical pendulum,  already studied by Huygens in the seventeenth   century,
is a famous  example of integrable system. 
In the  language of symplectic geometry it is described by  the 
cotangent bundle of $S^2$.  
Let $(\theta,\varphi)$ be the standard spherical angles, where $\varphi$ is
the rotation angle around the vertical axis and $\theta$
is the angle from the North Pole. Let $(\xi_\theta,\xi_\varphi)$
be the cotangent conjugate variables on ${\rm T}^*S^2$.
Then the cotangent bundle $({\rm T}^*S^2, \omega_{{\rm T}^*S^2})$ endowed
with the canonical cotangent bundle symplectic form
and the Hamiltonian function
$$
f_1(\underbrace{\theta,\varphi}_{\tiny \textup{sphere}},
\underbrace{\xi_\theta,\xi_{\varphi}}_{\textup{fiber}})=
\underbrace{\frac{1}{2}\left((\xi_\theta)^2 +
    \frac{(\xi_\varphi)^2}{\sin^2\theta}\right)}_{\textup{kinetic
    energy}}+ \underbrace{\cos\theta}_{\textup{potential}}.
$$
is an integrable system by
considering the vertical angular momentum $f_2(\theta,\varphi,\xi_\theta,\xi_{\varphi})= \xi_{\varphi}.$ 
So in our modern language, $F=(f_1,f_2) \colon {\rm T}^*S^2 \to \mathbb{R}^2$ is an
integrable system. Note that the Hamiltonian function $f_1$ is smooth on the cotangent bundle
${\rm T}^*S^2$ (the apparent singularity given by $1/\sin^2\theta$ is simply an artifact of
the choice of spherical coordinates).

\subsection{The topology of the fibers}

A point $m$ at which  $\mathcal{X}_{f_1}(m),\ldots,\mathcal{X}_{f_n}(m)$ of ${\rm T}_mM$
are linearly dependent is a \emph{singularity} of $F$. If $m$ is not
a singularity,  it is \emph{regular}. If a fiber $F^{-1}(c)$ contains some singularity we call it \emph{singular}. 
If  $F^{-1}(c)\neq \varnothing$ contains no singularity we call it \emph{regular}.

It follows from the definition of  integrable system, by simply following
the flows of the Hamiltonian vector fields of its components, that if $C$ is a connected component
of a regular fiber $F^{-1}(c)$, and if in addition 
the Hamiltonian vector fields $\mathcal{X}_{f_1},\ldots, \mathcal{X}_{f_n}$ 
are complete on $F^{-1}(c)$, then $C$ is diffeomorphic to $\mathbb{R}^{n-k} \times
\mathbb{T}^k$, where $\mathbb{T}^k:=(S^1)^k$ is a $k$\--dimensional torus.

\begin{figure}[h]
\centering 
\includegraphics[width=0.4\textwidth]{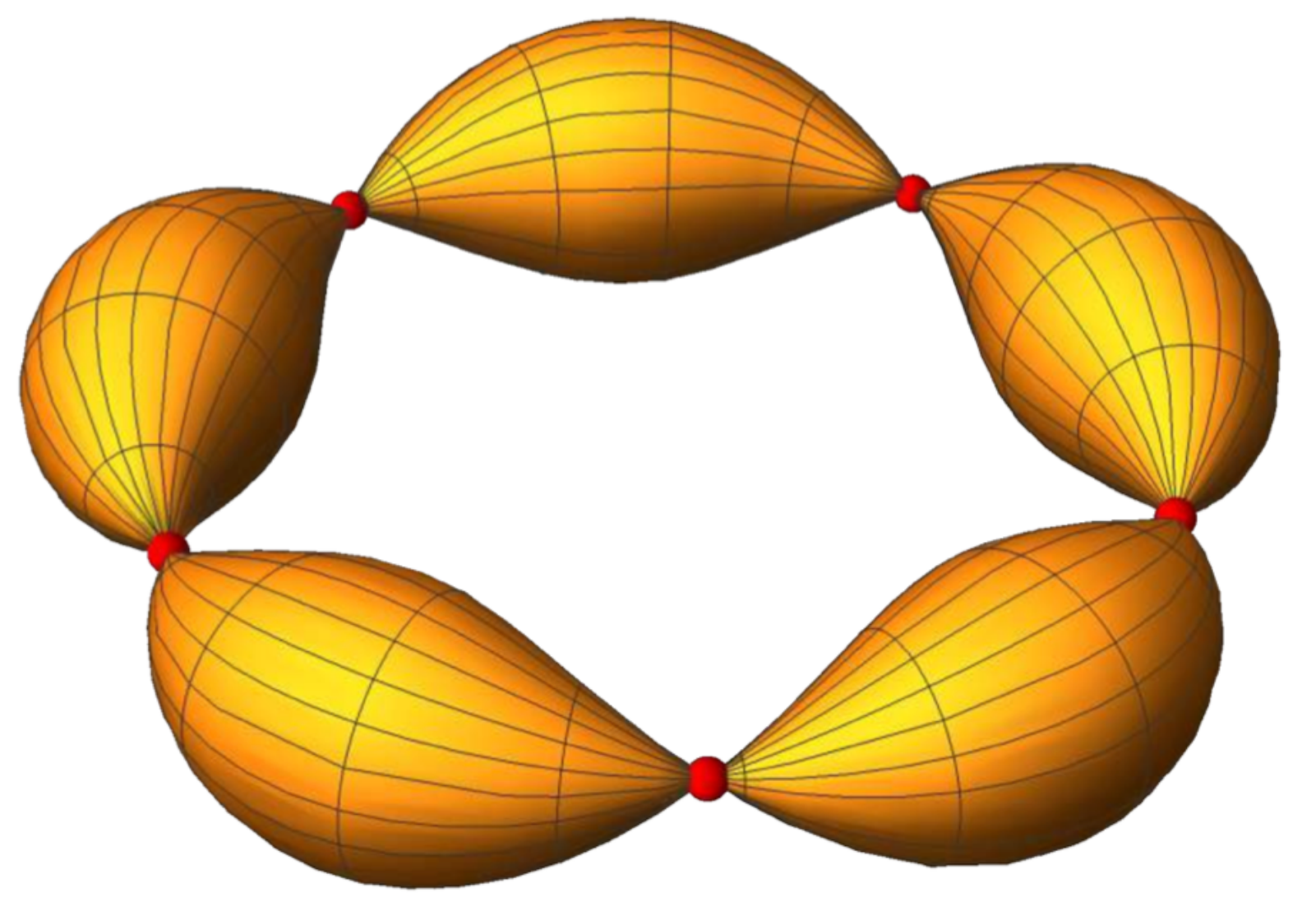}
 \caption{A wedge of $5$ spheres, or equivalently, a torus with $5$ pinches. Integrable systems
 $F \colon M \to \mathbb{R}^2$ can have singular fibers like these (generally, fibers with
 any number of pinched points exist). The symplectic geometry of this type of fibers has been
 understood only  recently, that is, that foliated neighborhoods
 of the fiber have been classified up to symplectomorphisms which preserve the foliation
 structure. These are typical singular fibers of non\--simple semitoric systems (in Section~\ref{nss}).}
\label{figure1}
\end{figure}

If the regular fiber $F^{-1}(c)$ is compact, then 
$\mathcal{X}_{f_1},\ldots, \mathcal{X}_{f_n}$  are complete on $F^{-1}(c)$, and $C$ is diffeomorphic to the $n$\--dimensional torus $\mathbb{T}^n$. If $F^{-1}(c)$ is both compact and connected, 
it is diffeomorphic to $\mathbb{T}^n$. Moreover,
$\omega$ vanishes along it, so $F^{-1}(c)$ is a Lagrangian submanifold 
(i.e. a submanifold on which $\omega$ vanishes); it has been traditionally called  a \emph{Liouville torus}. One often refers to $F \colon M \to \mathbb{R}^n$ as a singular Lagrangian fibration by Liouville tori, where ``singular" emphasizes 
that  $F$ may have singular fibers of various kinds:  for example  tori 
of dimension $m\in \{0,\ldots, n-1\}$, but also with more complicated topology, such as wedges of $2$\--spheres,
as in Figure~\ref{figure1}.

 A useful tool for proving results about integrable systems is
 Morse theory for the components $f_i$ of $F$. However, at least from the point of
 view of symplectic geometry, the applications of Morse theory are limited
 unless $M$ is compact or  $f_i$ is proper,  meaning that the preimages under $f_i$ of compact subsets of $\mathbb{R}$  are compact.    Furthermore, most of the 
  global \emph{symplectic} classifications of integrable systems known to date require that $F$
  is at least a proper map (which is automatic if at least one of the $f_i$ is proper). For 
   these reasons, unless otherwise stated, in this paper 
  \emph{integrable systems  $F$ are assumed to be proper maps into $\mathbb{R}^n$.} In particular
  this implies that all fibers of $F$ are compact.
   
  A different question is under which conditions the fibers of $F$
  are connected.  The question of whether $F^{-1}(c)$ 
  is connected can be formulated  as the question of whether the solution set to the system of equations
  $f_1=c_1,\ldots, f_n=c_n$ is or not connected, where $c=(c_1,\ldots,c_n)$. 
    Even if $M=\mathbb{R}^{2n}$ and the
  formulas defining the $f_i$ are simple, say quadratic formulas or polynomials, the question
  can be  challenging (it has  relations to real algebraic geometry).  At least in symplectic geometry, most approaches to this
  question have employed some version of Morse theory~\cite{atiyah,PeRaVNa,VN07}.

 \subsection{Action\--angle coordinates} \label{acc}
  
 Finite dimensional integrable systems have been studied from
 different view points and the literature, both in physics and mathematics, is extensive.
 Nonetheless, at least from the point of view of symplectic geometry, our  knowledge
 is  limited. 
 
 One of the few general theorems about the symplectic geometry 
 of integrable systems is the existence of action\--angle
 coordinates~\cite{arnold, mineur}  which says that each regular fiber, in addition to being
  diffeomorphic  to 
  $\mathbb{T}^n$ (recall that we are assuming throughout that $F$ is proper), sits inside of the cotangent
  bundle $T^*\mathbb{T}^n$ as the zero section, and in a neighborhood of the
fiber the integrable system  has the normal
form $F \colon T^*\mathbb{T}^n \to \mathbb{R}^n$, where
$$
F\underbrace{(x_1,\ldots, x_n}_{\textup{angle}},\underbrace{\xi_1,\ldots, \xi_n}_{\textup{action}})
=\underbrace{(\xi_1,\ldots,\xi_n)}_{\textup{action}}.
$$
To be more precise, we are assuming in this statement that we
have restricted ourselves to an adequate
invariant open set, so we can always assume that the fibers of $F$ are connected.  

 The existence of \emph{global} action\--angle coordinates was analyzed by
 Duistermaat in~\cite{Du}, a paper which may be considered to mark the beginning
 of the symplectic global theory of finite dimensional integrable systems. Since
 then this global problem has been analyzed in different contexts, see for
 example the case of non\--commutative integrable systems~\cite{RLF}.

\subsection{Linearization of non\--degenerate singularities} \label{s25}

In this paper we always assume that the singularities of $F \colon M \to\mathbb{R}^n$ are \emph{non\--degenerate}, we refer to~\cite{BAMS11} for the precise
definition. This notion is a vector\--valued extension of the condition of being Morse
non\--degenerate for real valued functions $M \to \mathbb{R}$ but it is 
more technical to describe.  The condition is satisfied by examples such
as the Jaynes\--Cummings model and the coupled angular momenta.

Under the condition of
non\--degeneracy, Eliasson described all local models for these singularities,
in what is considered one of the most influential results of the subject. His
result says that non\--degenerate singularities are 
linearizable~\cite{eliasson,eliasson-these, christophesan} 
in the sense that there exist coordinates  $(x_1,\ldots,x_n,\xi_1,\ldots,\xi_n)$ near
the singular point $m$ in which $m=(0,\ldots,0)$,  the symplectic form
$\omega$ has the expression $\omega= \sum_{i=1}^n{\rm d}x_i\wedge {\rm d}\xi_i,$ 
and there exist functions $q_1,\ldots, q_n$ of $(x_1,\ldots,x_n,\xi_1,\ldots, \xi_n)$  such that 
the integrable system $F=(f_1,\ldots, f_n)$ satisfies the Poisson bracket equation
$\{f_j,q_i\}=0,$ for all indices  $i,\,j$, where $q_i$ is one of
the following possibilities: 
\begin{enumerate}
\item
{elliptic type}: $q_i=\frac{x_i^2 + \xi_i^2}{2};$
\item
{hyperbolic type}: $q_i=x_i\xi_i$;
\item
{real type}: $q_i=\xi_i$;
\item
{focus\--focus type}: $q_i=x_i\xi_{i+1} - x_{i+1}\xi_i$ followed by
 $q_{i+1}=x_i\xi_i+x_{i+1}\xi_{i+1}$.
\end{enumerate}
If there are no components of hyperbolic type  the Poisson bracket equation  can be written 
as $$(F-F(m)) \circ \varphi=  g \circ (q_1,q_2,\ldots,q_n),$$ 
where
 $\varphi=(x_1,\ldots, x_n,\xi_1,\ldots,x_n,\xi_n)^{-1}$ and  $g$ is a diffeomorphism from
 a small neighborhood of $(0,\ldots,0)$ into another such neighborhood
 such that $g(0,\ldots,0)=(0,\ldots,0)$. For simplicity, usually we assume that $F(m)$ is the origin. 
 With few exceptions, in this paper we will restrict the discussion to
 singularities with no components of type (2). In this case Eliasson's theorem says that near any
 singularity there are symplectic coordinates in which $F$
 has the normal form $(q_1,\ldots, q_n)$ up to (translations and) composition
 by a local diffeomorphism.     One reason to rule out hyperbolic components  is that their appearance
 makes it difficult to construct global symplectic
 invariants.

\subsection{Singularities in dimension $4$} \label{d4}
Dimension $4$ is  the dimension (other than $2$) for which we have the best understanding
of the symplectic geometry of integrable systems.  The possibilities for $q_1$ and $q_2$ are:
 \begin{itemize}
 \item
$m$ is regular (rank $2$): $q_1=\xi_1$ and $q_2=\xi_2$;
\item
$m$ is transversally\--elliptic (rank $1$):   $q_1 = \frac{x_1^2 + \xi_1^2}{2}$ and $q_2 = \xi_2$; 
 \item
$m$ is {elliptic\--elliptic} (rank $0$): $ q_1 = \frac{x_1^2 + \xi_1^2}{2}$ and $q_2 = \frac{x_2^2 +
\xi_2^2}{2}$;
\item
$m$ is focus\--focus (rank $0$): $q_1=x_1\xi_2 - x_2\xi_1$ and $q_2 =x_1\xi_1+x_2\xi_2;$ 
 \item
$m$ is elliptic\--hyperbolic (rank $0$):  $q_1 = \frac{x_1^2 + \xi_1^2}{2}$  and $q_2 = x_2\xi_2$; 
 \item
$m$ is hyperbolic\--hyperbolic (rank $0$): $q_1=x_1\xi_1$ and $q_2=x_2\xi_2;$
\item
$m$ is transversally\--hyperbolic (rank $1$): $q_1=x_1\xi_1$ and $q_2=\xi_2$.
\end{itemize}
If we rule out hyperbolic components (last three possibilities), in coordinates
$(x_1,x_2,\xi_1,\xi_2)$ for which $\omega={\rm d}x_1\wedge {\rm d}\xi_1+{\rm d}x_2\wedge
{\rm d}\xi_2$, the  system $F$ has the 
form $(q_1,q_2)$ up to some local diffeomorphism. 
\subsection{Classifications and isomorphisms}

A leading goal of much of the recent research in symplectic geometry of
integrable systems has been to construct objects (numbers, functions, polytopes, 
etc) which are invariant by isomorphisms, in terms of which a class
of integrable systems can be classified up to these
isomorphisms. Examples of such classes could be: those which take place on
a $2$\--dimensional phase space, those whose
fibers are submanifolds, those for which the Hamiltonian vector field
associated to each component generates a periodic flow, etc.

 Here an \emph{isomorphism}
between integrable systems $F \colon (M,\omega) \to \mathbb{R}^n$ and
$F \colon (M',\omega') \to \mathbb{R}^n$ usually refers to a  diffeomorphism 
of phase space $f \colon M \to M'$ which
preserves the foliation structure by leaves induced
by the system viewed as a fibration (essentially meaning that $f^*F'= g \circ F$ for some smooth function $g$) 
and the symplectic structure (that is, $f^*\omega=\omega')$. 

This type of classification of integrable systems 
up to this notion of isomorphism is often referred to as a 
\emph{symplectic classification}, to emphasize the contrast with other
classifications in which one is interested in a notion of isomorphism which
does not necessarily have to preserve the symplectic structure. Such classifications
may be of a differentiable or topological nature instead.  There
are many works in these and related directions, 
see for instance~\cite{BS,BAMS11, Z} and the references therein.

At least from the point of view of spectral theory and quantization
of integrable systems, the most useful classifications need to be symplectic. 
 In dimension two 
there is a complete classification due to to Dufour\--Molino\--Toulet~\cite{DMT}.  
In dimensions higher than two  little is known, even less classifications, 
with two notable exceptions: toric integrable sytems on compact manifolds of
any dimension, and semitoric integrable systems (on compact or noncompact manifolds) but only 
in dimension $4$.  We will discuss these classifications in the upcoming sections.

Since symplectic classifications of integrable systems
have only been achieved in a few cases, often the specifics
of the notion of isomorphism have been tailored to each case, to obtain the most optimal form of
a classification. The notion above, depending on the context, may not be the most suitable.  
For example,  if we consider 
systems $F=(f_1,f_2) \colon M \to \mathbb{R}^2$ in which   the flow of $\mathcal{X}_{f_1}$ is periodic, the equation ``$f^*F'= g \circ F$" should be replaced by the more
specific equation $$f^*(f'_1,f'_2)=(f_1,h(f'_1,f'_2))$$ for some smooth function $h$ with 
$\frac{\partial h}{\partial f_2}>0$ as in~\cite{PeVN09}. For simplicity we will
not dwell on this; a discussion appears  in~\cite[Sections 5 and 6]{BAMS11}.  In Section~\ref{semi} we will
see a classification (up to this notion of isomorphism) 
of systems of this type, called \emph{semitoric}, under some  conditions.

\section{Toric integrable systems} \label{toricsection}

One of the branches of symplectic geometry which underwent a significant growth in the 1980s 
was the study  of Hamiltonian torus actions.  Effective Hamiltonian actions of tori of dimension $n$ on compact connected symplectic manifolds of dimension $2n$ were classified in~\cite{delzant}.  Such an action can be viewed as an integrable system  on a compact manifold for which all of its components generate periodic flows of the  same period, say $2\pi$; these systems are  called \emph{toric}. 

\subsection{The periodicity condition on the Hamiltonians}

We say that an integrable system 
 $$F=\underbrace{(f_1,\ldots,f_{n})}_{\textup{induce action of}\,\,\mathbb{T}^{n}} \colon M \to \mathbb{R}^n$$ 
on a symplectic $2n$\--dimensional manifold 
is \emph{toric} if  the Hamiltonian vector fields $\mathcal{X}_{f_1},\ldots, \mathcal{X}_{f_n}$
generate periodic flows of the same period, say $2\pi$, and  the action of $\mathbb{T}^n$ on $M$ produced
by concatenating these flows is effective.

Unless otherwise stated, in this paper 
we will only consider \emph{toric integrable systems on compact connected manifolds}. Indeed compactness is
a crucial condition  for the results in this section and for the applications in Section~\ref{spec1}.

The periodicity implies that the singularities of toric integrable systems 
cannot have focus\--focus or hyperbolic type components, that is, if 
$m=(0,\ldots,0)$ and $\omega=\sum_{i=1}^n{\rm d}x_i \wedge {\rm d}\xi_i$,
locally in a neighborhood of $m$ the integrable system must have the form
$$
F(x_1,\ldots,x_1, \xi_1, \ldots, \xi_n)=\Big(\underbrace{\frac{x_1^2 + \xi_1^2}{2},\ldots,\frac{x_k^2 + \xi_k^2}{2}}_{\textup{elliptic type}}, \overbrace{\xi_{k+1},\ldots,\xi_{n}}^{\textup{real type}}\Big). 
$$

Toric integrable systems have connected fibers, a fact  known
as Atiyah's connectivity~\cite{atiyah} in the more general 
context of Hamiltonian $\mathbb{T}^m$\--actions, $m \in \{1,\ldots,n\}$, 
on compact connected symplectic $2n$\--dimensional manifolds (there
are  also extensions of this to certain infinite dimensional settings, see for instance~\cite{at-pr, HHJM}). This fact is
closely  related to the classification of toric integrable systems which we discuss in Section~\ref{delzant}.

All fibers of $F$ are diffeomorphic to tori of varying dimensions 
$\mathbb{T}^k$,  $k \in \{0,\ldots, n\}$, which is not the case for general integrable systems as we will
see in a moment (earlier we say a typical singular fiber in Figure~\ref{figure1}). For example, in dimension $4$ the possibilities are:
\begin{eqnarray}
&&F(x_1,x_2,\xi_1,\xi_2)=\Big(\frac{x_1^2 + \xi_1^2}{2},\frac{x_2^2 + \xi_2^2}{2}\Big) \Longrightarrow   F^{-1}(F(m))=\underbrace{\{m\}}_{\textup{elliptic\--elliptic singularity}};    \nonumber \\
&&F(x_1,x_2,\xi_1,\xi_2)=\Big(\frac{x_1^2 + \xi_1^2}{2},\xi_2\Big)  \Longrightarrow 
\underbrace{F^{-1}(F(m))\simeq S^1}_{\textup{transversally\--elliptic singularities}}; \nonumber \\
&&F(x_1,x_2,\xi_1,\xi_2)=(\xi_1,\xi_2)  \Longrightarrow  
\underbrace{F^{-1}(F(m)) \simeq \mathbb{T}^2}_{\textup{regular points}}. \nonumber
\end{eqnarray}
The fibers corresponding to these local models are illustrated in Figure~\ref{figure2}, in terms of the
image under $F$ of the singularity.

\subsection{Example of a toric integrable system: rotation on complex projective space}

Consider on $S^2$ the standard area form  $\omega={\rm d}\theta \wedge {\rm d}h$, where $\theta$
is the angle and $h$ is the height of a point in $S^2$.  The height function $F(\theta,h)=h$ is a Hamiltonian on $S^2$, which defines a toric
integrable system. This system corresponds to 
the rotational $S^1$\--action on $S^2$ about the $z$\--axis, which has momentum
map precisely equal to $F$. Clearly $F(S^2)=[-1,1]$.  

As a generalization  to higher dimensions, 
consider the complex projective space 
$\mathbb{C}P^n$ endowed with a $\lambda$\--multiple ($\lambda>0$) of the Fubini\--Study form. 
Then $F \colon \mathbb{CP}^n \to \mathbb{R}^n$   given by
$$
 F([z_0:z_1:\ldots:z_n])=
\Big(\frac{\lambda |z_1|^2}{\sum_{i=0}^n|z_i|^2},\ldots, \frac{\lambda |z_n|^2}{\sum_{i=0}^n|z_i|^2}\Big).
$$
is a toric integrable system, induced by the Hamiltonian $\mathbb{T}^n$\--action  by rotations defined 
by the formula
$({\rm e}^{{\rm i}\theta_1}, \ldots 
{\rm e}^{{\rm i}\theta_n})[z_0:z_1,\ldots:z_n]=
[z_0: {\rm e}^{{\rm i}\theta_1}z_1: \ldots {\rm e}^{{\rm i}\theta_n}z_n].
$
 If ${\rm e}_1=(1,0,\ldots,0),\ldots,{\rm e}_n=(0,\ldots,0,1)$, are
the standard basis vectors in $\mathbb{R}^n$, then the image of the integrable system is the set
$
F(\mathbb{C}P^n)=\textup{convex hull }\{0, \lambda {\rm e}_1,\ldots,\lambda {\rm e}_n\}.
$
Notice that $F$ has an interesting property: its image, which is a subset of $\mathbb{R}^n$,
is a convex polytope. Moreover, this convex polytope is the convex hull of the fixed point
set of the $\mathbb{T}^n$\--action on $\mathbb{C}P^n$. This property turns out to be a general
property of toric integrable systems, as we explain next.

\subsection{Classification} \label{delzant}

One of the fundamental theorems of equivariant symplectic geometry,  due to Atiyah~\cite{atiyah} and
Guillemin\--Sternberg~\cite{gs} says that, if $M$ is compact and connected, the image $F(M)$ is a convex
polytope in $\mathbb{R}^n$, in fact, obtained by a recipe: it is the convex hull of
the images of fixed points of the Hamiltonian action of the $n$\--torus on $M$ induced by concatenating 
the flows of the $f_i$. This polytope has the property that if two 
toric integrable systems are isomorphic then their associated images coincide 
(up to translations and composition with a matrix in ${\rm GL}(n,\mathbb{Z})$).

 In fact, the Atiyah\--Guillemin\--Sternberg convexity 
theorem\footnote{Convexity and related 
properties of Hamiltonian actions were later studied in different contexts,
for example see for instance~\cite{Alek, FlRa, Ki, OrRa, We2}, to name a few.}, as  is usually
referred this result, applies to general Hamiltonian $m$\--dimensional torus actions, for any 
$m\in \{1, \ldots, n\}$, where $2n$ is the dimension of $M$. In this case 
$F$ would be the action momentum map and has  $m$ components. 
Of course only in the case when $m=n$  this represents simultaneously a Hamiltonian torus action
and an integrable system on $M$. It is in this case when their result can be used as a stepping stone for obtaining a  classification of toric integrable
systems on compact manifolds.

Indeed,  
shortly after this result, Delzant showed that the polytopes obtained as images of toric integrable 
systems are of a special type: simple, rational, and smooth. Now these polytopes
 are  called \emph{Delzant polytopes}; the
essential condition for a Delzant polytope in $\mathbb{R}^n$ is that there are precisely $n$ edges
meeting at each vertex and the normal vectors to the facets meeting at the vertex form
a basis of the integral lattice.

\begin{figure}[h]
\centering 
\includegraphics[width=0.35\textwidth]{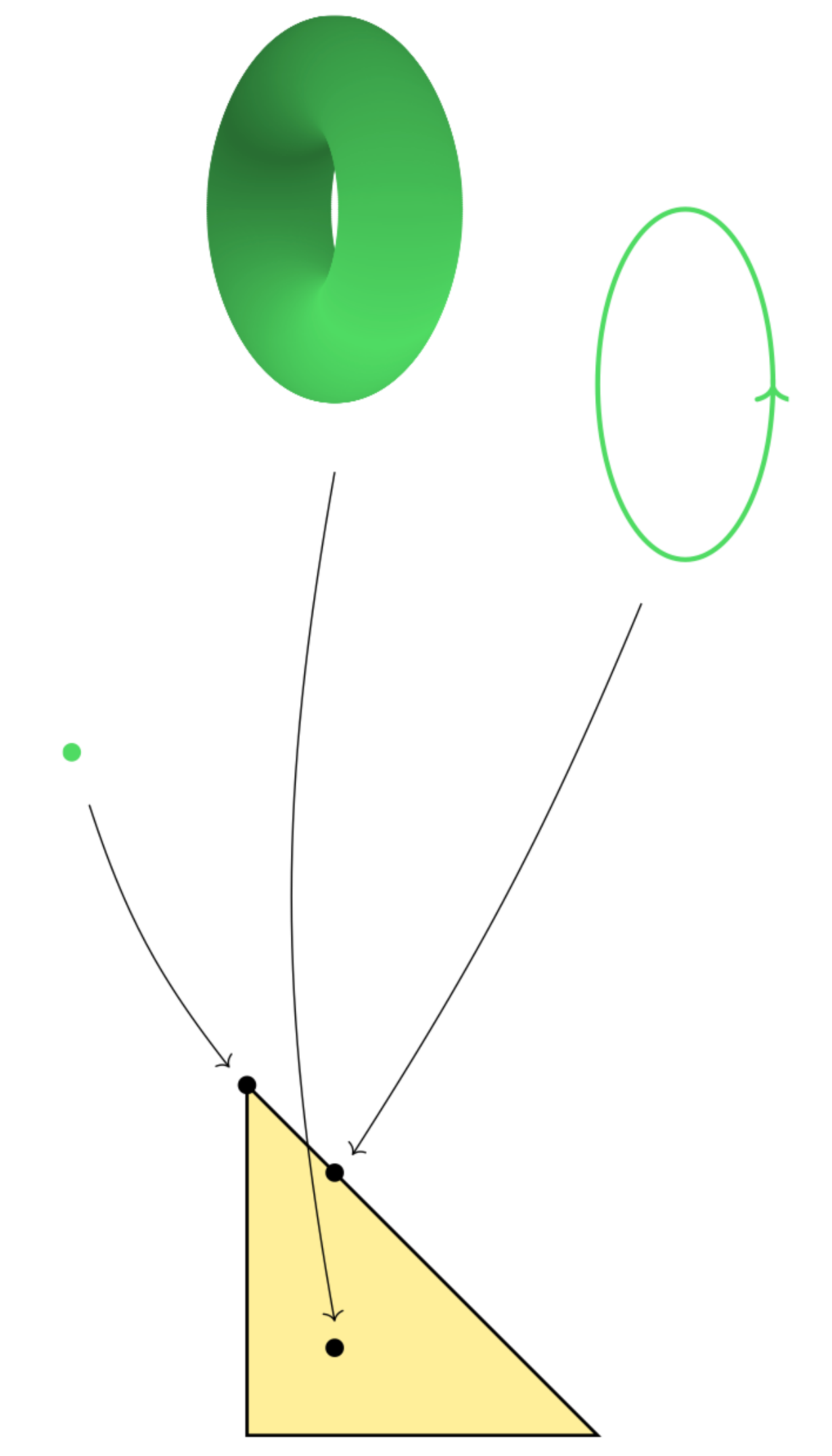}
 \caption{The information of a toric integrable system $F \colon (M,\omega) \to \mathbb{R}^n$ on a compact
 connected manifold can be read off from the polytope $F(M)=\Delta$. In particular the topology of a fiber over a point $p \in \Delta$ is  encoded by the dimension of the face of $\Delta$ containing $p$. The fiber over a vertex
 of $\Delta$ is a point, while the fiber over a point in an edge of $\Delta$ is a circle. If $p$ is any point
 in the interior of the polygon, the fiber is a torus.}
\label{figure2}
\end{figure}

Moreover, Delzant proved that such polytopes are in bijective correspondence with toric integrable systems on compact manifolds:
$$
 \textup{toric system} \underbrace{\rightsquigarrow}_{\textup{information lost?}}  \overbrace{\textup{its image} \,\,\Delta}^{\textup{Delzant polytope}}  \underbrace{\rightsquigarrow}_{\textup{Delzant's Theorem}}  \textup{toric system}
$$

Here
by ``bijective correspondence" we mean:
\begin{itemize}
\item
 the \emph{uniqueness} statement that two  systems 
 $F \colon (M,\omega) \to \mathbb{R}^n$ and $F' \colon (M',\omega')  \to \mathbb{R}^n$
 are isomorphic
if and only if they have the same convex polytope as image (up to translations and ${\rm GL}(n,\mathbb{Z})$ transformations), so the answer to the question in the first arrow of the diagram is \emph{no} (up to isomorphisms), and,
\item
the \emph{existence} statement which says that starting from any Delzant polytope $\Delta$ in $\mathbb{R}^n$
one can construct a toric integrable system $F \colon (M,\omega) \to
\mathbb{R}^n$  on a compact connected symplectic $2n$\--dimensional 
manifold whose image is   $\Delta$. This construction, which corresponds
to the second arrow of the diagram,
 can be achieved via the method of \emph{symplectic reduction}. 
\end{itemize}

Since $\Delta$ classifies $F \colon (M,\omega) \to \mathbb{R}^n$, one can know everything
about $F$, up to isomorphisms, from $\Delta$.  In particular 
the fiber structrure of $F$ can be read off from the polytope $\Delta$:  the fiber of $F$ over $p \in \Delta$
is diffeomorphic to an $k$\--dimensional torus, where $k$ is the dimension of the lowest
dimensional face of $\Delta$ such that $p \in \Delta$. For example, in the case of $\mathbb{CP}^2$
discussed earlier,  $\Delta$ is a $2$\--dimensional simplex with boundary a triangle. The fibers
of $F$ over the vertices of the triangle are points (elliptic\--elliptic singularities), over the edges of the triangle are diffeomorphic to $S^1$ (transversally\--elliptic singularities),
 and over any point in the interior of the triangle they are $2$\--tori (regular points); this is depicted
 in Figure~\ref{figure2}.

Delzant's classification has  been 
 a precursor for the study of more general integrable systems. 
 Since then, there have been  classifications of
systems of toric type in a variety contexts. For instance, on noncompact
manifolds~\cite{KL2015}, and on log\--symplectic manifolds~\cite{GMP2014,LGPR2017}.  Delzant's
result has been extended  to general symplectic torus 
actions in some cases~\cite{benoist, DP1,Pe2}.
 There is also a natural connection between the study of toric integrable systems and
 toric varieties in algebraic geometry, see~\cite{duispel3} for a detailed study.

 \section{Semitoric integrable systems} \label{semi}
 
 Semitoric integrable systems form a class of integrable systems which generalizes the class
 of toric integrable systems. They have been the focus of intense research activity in the past ten years or so. 
 They were classified in dimension $4$ under some conditions in~\cite{PeVN09,PeVN11}.  
 The main property of a semitoric system is that all the functions that define it, but one, 
 generate periodic flows. So  toric systems are semitoric. At least from the 
 point of view of symplectic geometry, the main difference with toric systems is that semitoric
 systems can have focus\--focus singularities and the fibers containing them are not diffeomorphic to tori, 
for instance as  in Figure~\ref{figure1}   (while the fibers of toric  
 systems are diffeomorphic to tori of varying dimensions, as  in 
 Figure~\ref{figure2}). As we will see, the appearance of
 this type of singular fiber makes the symplectic geometry of these systems very rich.
  
 The Jaynes\--Cummings model is a well known semitoric system. Its fiber structure looks like that of a toric system, with the exception that
 it has one fiber containing a focus\--focus singularity, which is a torus pinched at exactly one point. The coupled angular momenta
 is another example.

\subsection{The periodicity condition on all but one Hamiltonian}

 An integrable system $$F=\underbrace{(f_1,\ldots,f_{n-1}}_{\textup{induce action of}\,\,\mathbb{T}^{n-1}},f_n) \colon M \to \mathbb{R}^n$$ 
 is \emph{semitoric} if   the Hamiltonian
 vector fields $\mathcal{X}_{f_1},\ldots, \mathcal{X}_{f_{n-1}}$ generate periodic flows of the
 same period, say $2\pi$, and the action of $\mathbb{T}^{n-1}$ on $M$ produced by concatenating these flows
 is effective.  As  earlier,  we also require that the singularities of $F$ are non\--degenerate 
 and do not have hyperbolic components.  In the current theory there is also
 the assumption that $f_1,\ldots, f_{n-1}$ are all proper.

 If $n=2$ these conditions  mean (by Section~\ref{d4}) that $f_1 \colon M \to \mathbb{R}$ 
 is a proper momentum map for an effective Hamiltonian $S^1$\--action, 
  and the local models of $F=(f_1,f_2)$ near a point $m$ are, in some
 coordinates $(x_1,x_2,\xi_1,\xi_2)$ in which $m=(0,0,0,0)$ and
 $\omega={\rm d}x_1\wedge {\rm d}\xi_1+{\rm d}x_2\wedge {\rm d}\xi_2$,
 one of the following:
 \begin{eqnarray}
&&F(x_1,x_2,\xi_1,\xi_2)=\Big(\frac{x_1^2 + \xi_1^2}{2},\frac{x_2^2 + \xi_2^2}{2}\Big) \Longrightarrow   F^{-1}(F(m))=\underbrace{\{m\}}_{\textup{elliptic\--elliptic singularity}};    \nonumber \\
&&F(x_1,x_2,\xi_1,\xi_2)=\Big(\frac{x_1^2 + \xi_1^2}{2},\xi_2\Big)  \Longrightarrow 
\underbrace{F^{-1}(F(m))\simeq S^1}_{\textup{transversally\--elliptic singularities}}; \nonumber \\
&&F(x_1,x_2,\xi_1,\xi_2)=(\xi_1,\xi_2)  \Longrightarrow  
\underbrace{F^{-1}(F(m)) \simeq \mathbb{T}^2}_{\textup{regular points}}; \nonumber \\
&&F(x_1,x_2,\xi_1,\xi_2)= (x_1\xi_2 - x_2\xi_1, x_1\xi_1+x_2\xi_2)  \Longrightarrow  
\underbrace{F^{-1}(F(m)) \simeq \textup{pinched}\,\, \mathbb{T}^2}_{\textup{a few focus\--focus singularities inside}}. \nonumber
\end{eqnarray}
The singular fibers, which are connected in this case and also under
 slightly weaker conditions~\cite{PeRaVNa, VN07}, are either points, circles, or tori pinched
 at a finite amount of points (i.e. wedges of $2$\--spheres as in Figure~\ref{figure1}). This last type of fiber does not appear
 in toric integrable systems.
 
  If $n=2$ we have a complete understanding of the global symplectic geometry of semitoric  systems
  under the assumptions above: they were classified about ten years ago 
 in~\cite{PeVN09,PeVN11} under the extra requirement that  each fiber of $f_1$, and hence of $F=(f_1,f_2)$, contains at most
one focus-focus singularity; this is often called \emph{simplicity}. So each singular fiber of $F$ could be pinched at most once.
 This requirement was recently removed in~\cite{PPT19}, hence in particular allowing fibers of the
 form shown in Figure~\ref{figure1}. We will discuss this recent extension in Section~\ref{nss}, but before that section we assume that all systems satisfy this simplicity condition. 
 
\subsection{Example of a simple semitoric integrable system: the Jaynes\--Cummings model}

A well known example of semitoric system, which is also one of the few for
which the symplectic invariants  have been explicitly computed, 
is the Jaynes\--Cummings model~\cite{JC63, C65} from
physics, also referred to  as the coupled spin\--oscillator. This 
is a system which models simple physical phenomena, and is given by 
$$F(x,y,z,u,v):=\Big(\underbrace{f_1= \frac{u^2+v^2}{2} + z}_{\textup{periodic flow}}, f_2=\frac{ux+vy}{2}\Big)$$ in coordinates 
$(x,y,z,u,v) \in S^2 \times \mathbb{R}^2$, where $S^2$ inherits the coordinates from
the usual inclusion $S^2 \subset \mathbb{R}^3$. Note that $f_1$ corresponds to the momentum 
map for the Hamiltonian $S^1$\--action by rotations on the plane about the origin.
All the singularities of
$F$ with the exception of $m=(0,0,1,0,0)$ are elliptic or transversally
elliptic, while $m$ itself is a focus\--focus singularity~\cite{PeVN12}. Hence the singular
fiber containing $m$ is  a torus pinched at $m$. 

There are other famous systems which are semitoric, notably the coupled
angular\--momenta. An example which is not semitoric strictly
speaking, but for which one component still generates a periodic flow is the
spherical pendulum in Section~\ref{nm}. In this case 
the system fails to satisfy our definition of semitoric system because  the Hamiltonian
generating a periodic flow is not proper, even though the joint map $F$ is. 

 While there is no general symplectic 
classification for examples like the spherical pendulum in which the periodic
component is not proper, some first steps were taken in~\cite{PeRaVNa, PeRaVNb} to 
construct a polygonal invariant out of the image of $F(M)$ which resembled
the convex polygon of Atiyah\--Guillemin\--Sternberg~\cite{atiyah, gs} and its
later generalization to semitoric systems~\cite{PeVN09, VN07}.

\subsection{Symplectic invariants and uniqueness of simple semitoric systems in dimension four} \label{uni}

\begin{figure}[h]
\centering 
\includegraphics[width=0.4\textwidth]{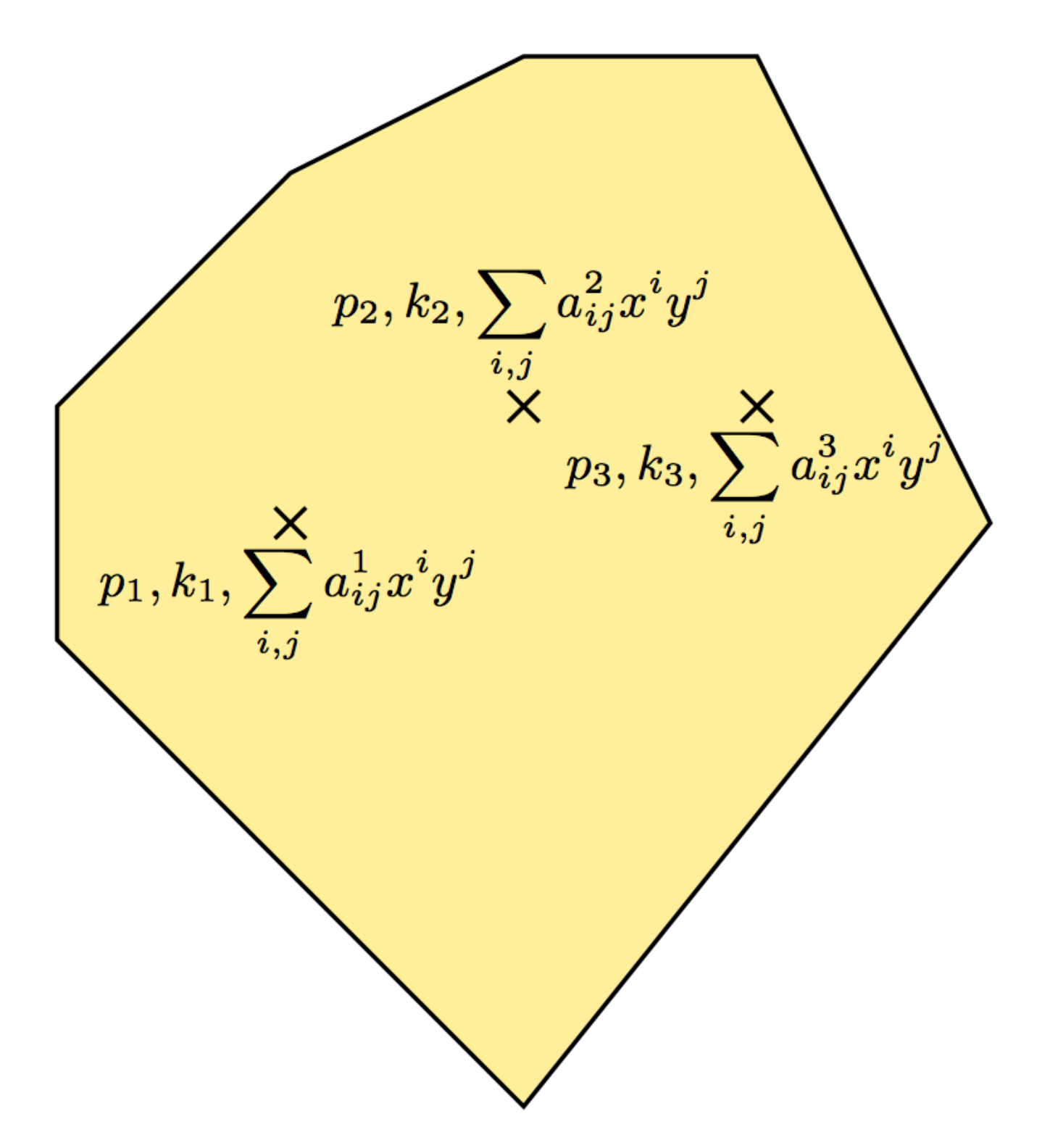}
 \caption{The figure shows a polygon $\Delta$ in the plane, with several points inside of it. Each point
 comes endowed with a label consisting of an integer and a formal Taylor series on variables $x,y$. This
 is the generic picture which captures the symplectic invariants of a simple semitoric integrable system on a $4$\--manifold. In this sense, the figure is a \emph{complete invariant}. That is, all of the symplectic geometry of the simple semitoric system $F \colon (M,\omega) \to \mathbb{R}^2$ can be read off from this figure. The uniqueness theorem for semitoric systems explained in 
  Section~\ref{uni} says that any two simple semitoric systems with the same complete invariant must be isomorphic.  However, this result does not tell us which polygons can appear, and whether there are
  restrictions on the positions of the points, the coefficients of the Taylor series, or the integers $k$. This is
  a separate, existence problem,  discussed in Section~\ref{existence}.}
\label{figure3}
\end{figure}

In dimension four, a simple semitoric system $F=(f_1,f_2)$ is determined~\cite[Theorem 6.2]{PeVN09} up to isomorphisms by a convex  polygon $\Delta$ with marked points $p_1,\ldots,p_n \in \Delta$,  each of which
comes with a label assigned  $$p_{\ell} \rightsquigarrow  \Big(k, \sum_{i,j} a_{ij}x^iy^j\Big) \in \mathbb{Z} \times \mathbb{R}[x,y],$$ where $k$ is an integer and  $\sum_{i,j} a_{ij}x^iy^j$
is a formal Taylor series on two variables.  These invariants are depicted in Figure~\ref{figure3}.

The main tools for proving this
 are the existence of action\--angle coordinates~\cite{arnold,Du} and the results on
linearization of non\--degenerate singularities in~\cite{eliasson, eliasson-these,MZ, christophesan}.  
Note that the $k$ and the coefficients $a_{ij}$ 
of the Taylor series will vary, in general, for each marked point. Notice that this is a \emph{uniqueness result},
while it does give no information on how many simple semitoric systems there are or how to construct them,
which would be an~\emph{existence result}. We well discuss this existence question shortly (in Section~\ref{existence}).

So the list of symplectic invariants of simple semitoric systems
consists of a polygon, decorated by $n$ points and $n$ 
labels\footnote{The number $n$ of marked points and their exact positions
inside of $\Delta$ (more precisely along the vertical line which contains them) are both symplectic invariants. With this point of view, one could say that there
are five invariants: the polygon $\Delta$ itself, the number of points inside of it, the position of 
each point inside of $\Delta$,  the $n$-tuple of integer labels $k$, and the $n$\--tuple of
Taylor series, one per point. These were the five invariants introduced in~\cite{PeVN09} in terms
of which semitoric systems are classified in~\cite{PeVN09,PeVN11}; our
formulation in the present paper is equivalent.}. The marked points here correspond
to the images of the focus\--focus singularities of $F$ (there are finitely many
of them as shown in~\cite{VN07}).  Hence if the system is toric, then there
are no marked points and hence no labels, so we are left with a convex
polygon, and if $M$ is compact this recovers the classification result
of Delzant in Section~\ref{delzant}. 
Next we briefly explain the meaning of these invariants:

\begin{itemize}
\item
The convex polygon $\Delta$ 
is equal to $F(M)$ if $M$ is compact and the system is toric. In general,
$\Delta$ does not equal $F(M)$, instead it is obtained from $F(M)$ by making
vertical cuts along the focus\--focus values~\cite{VN07}, hence unfolding the singular affine
structure induced by $F$. These vertical
cuts can be made either upwards, or downwards, resulting in a family of polygons and
the actual polygon invariant is an orbit of such polygons by the action which takes into consideration
the cut directions; it is called the \emph{semitoric polygon invariant}, the precise notion is given
in~\cite[Section 4.3]{PeVN09}. In fact, $\Delta$
 satisfies some specific rationality properties as explained in~\cite[Section 4.1]{PeVN11}, which
generalize those of the Delzant polygons in the classification of toric integrable systems on compact
connected manifolds.
\item
The Taylor series $\sum_{i,j} a_{ij}x^iy^j$ encodes the  dynamical behavior of the
Hamiltonian vector fields  $\mathcal{X}_{f_1}$ and $\mathcal{X}_{f_2}$
as they approach the focus\--focus singularities of $F$, and symplectically
determines the semiglobal normal form of the fibration by Liouville tori near the singular fiber containing
the focus\--focus singularity corresponding to the point in the polygon with the label 
$(k, \sum_{i,j} a_{ij}x^iy^j)$. 
\item
The integer $k$, called the \emph{twisting index}, 
encodes the topology of  $F$
viewed as a singular Lagrangian fibration near the singular fiber 
containing the focus\--focus singularity corresponding to the point in the polygon with the label 
$(k, \sum_{i,j} a_{ij}x^iy^j)$.   This index can be thought of as encoding
 the way in which the semiglobal normal form determined by $\sum_{i,j}a_{ij}x^iy^j$ is glued
 into the integrable system, relative to the global integral affine structure induced by $F$. 
 The construction of $k$ is subtle and beyond
the scope of this paper,  we refer to~\cite[Section 5.2]{PeVN09} for those interested in
learning about it. The actual \emph{twisting index invariant} is given as an orbit
of the tuple of all twisting indices, say $(k_1,\ldots, k_n) \in \mathbb{Z}^n$, corresponding to the points 
$p_1,\ldots,p_n \in \Delta$ under some natural group action on them. 
\end{itemize}

To get a better idea of the result, we explain how to construct
$\sum_{i,j} a_{ij}x^iy^j$. Let $m$ be a focus\--focus
singularity of $F \colon M \to \mathbb{R}^2$. 
By Eliasson's theorem there are coordinates $(x_1,x_2,\xi_1,\xi_2)$ near $m$ 
in which $m=(0,0,0,0)$,  $\omega={\rm d}x_1 \wedge {\rm d}\xi_1 +{\rm d}x_2\wedge {\rm d}\xi_2$ 
and a local diffeomorphism $h$ such that 
$
h \circ {F}=(x_1\xi_2 - x_2\xi_1, x_1\xi_1  +x_2\xi_2).
$
This normalization extends to a neighborhood of $F^{-1}(F(m))$, we call it
$
(H_1,H_2).
$
Now we are going to construct two period maps. For each regular value  $c$ close to $F(m)$
choose  $A \in F^{-1}(c)$ and define
$\tau_2(c)$  to be the time it takes the Hamiltonian flow associated
to $\mathcal{X}_{H_2}$ leaving from $A$ to meet the Hamiltonian flow
associated to $H_1$ which passes through $A$, and let 
$\tau_1(c) \in \mathbb{R}/2\pi \mathbb{Z}$  be
the time that it takes to go from this intersection point back to $A$ along the Hamiltonian flow
line of $\mathcal{X}_{H_1}$, in this way closing the trajectory. 
The behavior of $\mathcal{X}_{H_2}$
 field as we approach $c$ becomes singular. We remove the singular
 behavior from the period maps, which is of a logarithmic nature, by defining
   $\sigma_1(c)=\tau_1(c)-{\rm Im}(\log(c))$ and 
 $\sigma_2(c)=\tau_2(c)+{\rm Re}(\log(c))$, which are well\--defined and
smooth, and from them  we construct the $1$\--form
 $
 \sigma_1{\rm d}c_1+\sigma_1{\rm d}c_2
 $
 which is closed near $m$ and hence exact: $\sigma_1{\rm d}c_1+\sigma_2{\rm d}c_2={\rm d} S$. The invariant $\sum a_{ij} x^i y^j$ is the Taylor series of $S$ (up to a $(\mathbb{Z}_2 \times \mathbb{Z}_2)$\--action, as 
 in~\cite{SeVN}). 
 
 The invariant $\sum a_{ij} x^i y^j$  is a semiglobal invariant~\cite{VN03}, 
 that is, it determines up to isomorphisms
 a neighborhood of $F^{-1}(F(m))$, meaning two integrable systems are isomorphic
 in a neighborhood of a fiber with a focus\--focus singularity if and only if they have the same Taylor series
associated to them. In fact,
 essentially any Taylor series of two variables (with a small restriction on a coefficient)
 can appear as a symplectic invariant, and hence $\sum a_{ij} x^i y^j$
 provides a symplectic classification a neighborhood of  the fiber containing $m$.

 In the case that there are $\lambda>1$ focus\--focus singularities in the same singular fiber, it has
 been recently shown~\cite{PeTa} that the symplectic invariant of a neighborhood of the singular
fiber containing the $\lambda$ focus\--focus singularities consists of $\lambda$ Taylor series:
$
 \sum_{i,j} a^1_{ij}x^iy^j,\ldots,  \sum_{i,j} a^{\lambda}_{ij}x^iy^j,
 $
which again encode the singular dynamics of the Hamiltonian vector fields
 $\mathcal{X}_{f_1}$ and $\mathcal{X}_{f_2}$ as they approach these $\lambda$ singularities.

\subsection{Calculations of symplectic invariants of certain semitoric systems}

The Taylor series invariant has been calculated in some examples; for the case of the Jaynes\--Cummings
system the  Taylor 
series at the only focus\--focus singularity  was calculated in~\cite{PeVN12} to be
$$\sum_{i,j} a_{ij} x^iy^j=\frac{\pi}{2}  x+ 5\log 2 y  +\mathcal{O}(2).$$

A detailed study of the
invariants of this system was later given in~\cite{ADH19}, see~\cite{AlHo} for a 
survey of results in this direction.

Let's briefly discuss another example: the case of the coupled angular momenta~\cite{SZ}. 
 Let $R_2 > R_1 > 0$. On $M=S^2 \times S^2$ with
coordinates $(x_1,y_1,z_1,x_2,y_2,z_2)$ we consider:
\begin{equation} \begin{cases} 
f_1 := R_1 z_1 + R_2 z_2 \\ \\ f_{2,t} := (1-t) z_1 + t (x_1 x_2 + y_1 y_2 + z_1 z_2) 
\end{cases}
\forall t \in [0,1] \nonumber
\end{equation}
with symplectic form $ -(R_1 \omega_{S^2} \oplus R_2 \omega_{S^2})$, where $\omega_{S^2}$
is the standard symplectic form on $S^2$. In this case there are two values $t^{-}<\frac{1}{2}<t^+$
such that  $(0,0,1,0,0,-1)$  is a singularity of 
$$
F_t:=(f_1,f_{2,t}) \colon S^2 \times S^2 \to \mathbb{R}^2
$$
which is elliptic\--elliptic 
between $0$ and $t^{-}$, degenerate at $t^{-}$, focus\--focus between $t^{-}$ and $t^+$, degenerate
at $t^+$, and elliptic\--elliptic again from $t^+$ to $1$. 

The integrable system $F_t$ is 
toric up to local diffeomorphisms (this is often called being of \emph{toric type}) for $t<t^{-}$
and $t>t^{+}$, and semitoric
for values of $t$ in between $t^{-}$ and $t^+$.  In~\cite{LFPe19} the polygon invariant
was computed for $t = 1/2$, and it was also shown that  there is a focus\--focus
singularity with label $k=0$, and if  $R_1 = 1, R_2 = 5/2$ the Taylor series is 
$$
\sum_{i,j} a_{ij} x^iy^j=\arctan\Big(\frac{9}{13}\Big) x + \left(\frac{7}{2} \log 2 + 3 \log 3 - \frac{3}{2} \log 5 \right) y + \mathcal{O}(2).
$$

For a complete study of the symplectic invariants of the coupled angular momenta see~\cite{ADH19b}. Further generalizations and other families of deformations were studied in~\cite{HP, LFPa}.

\subsection{The construction (i.e. existence) of  simple semitoric systems in dimension $4$} \label{existence}

It was shown in~\cite{PeVN11} that for each semitoric polygon $\Delta$ 
with interior marked points  $p_1,\ldots,p_n$, each labelled by a pair
$
(k, \sum_{i,j} a_{ij}x^iy^j),
$ 
  there exists a semitoric integrable system 
$F \colon M \to \mathbb{R}^2$ 
which has $\Delta$ and $p_1,\ldots, p_n$, each with a label $(k, \sum_{i,j} a_{ij}x^iy^j)$, 
as its symplectic invariants, as defined in Section~\ref{uni} (this is depicted in 
Figure~\ref{figure3}). Together with the uniqueness 
result in Section~\ref{uni}, this gives a symplectic classification of simple semitoric integrable systems. 
(The precise formulation of the classification in~\cite[Theorem 4.7]{PeVN11} involves considering 
certain natural group actions on the space of marked semitoric polygons, to avoid redundancies and making
sure the invariants are defined without ambiguity). We may illustrate this with the following diagram, while
keeping in mind that the first arrow was discussed in Section~\ref{uni}, while the second arrow is
the purpose of this section:
$$
\textup{simple semitoric syst.}  \underbrace{\rightsquigarrow}_{\textup{lost?}} 
  \overbrace{\Delta,p_1,k_1, \sum_{i,j} a^1_{ij}x^iy^j,\ldots, p_n, k_n, \sum_{i,j} a^n_{ij}x^iy^j}^{\textup{symplectic invariants}}  \underbrace{\rightsquigarrow}_{\textup{\S\ref{existence}}}  \textup{simple semitoric syst.}
$$

By Section~\ref{uni} the answer to the question in the first arrow of the diagram, as to whether
any data is lost in the process of finding the invariants,  is \emph{no} (up to isomorphisms). 
For brevity and to keep the discussion  nontechnical we have omitted 
what type of symplectic invariants can appear as an abstract list
of objects defined without reference to symplectic manifolds or integrable systems. In the toric
case, recall this list consisted of all Delzant polytopes in $\mathbb{R}^n$; in the semitoric
case, one needs to allow more flexibility on the polytopes that can appear, as well as have other
conditions imposed on the abstract spaces where the remaining symplectic invariants are: space
of formal Taylor series, $\mathbb{Z}$, etc.

The idea of proof of this existence result, which shows an explicit construction, is as follows. 
Basically one needs to glue the ``regular" part of the system (given by action\--angle coordinates)
with the ``singular" part, near the focus\--focus fibers (which involves Eliasson's coordinates,
as employed for the construction of the Taylor series invariant in Section~\ref{uni}):
\begin{itemize}
\item
In the first step $\Delta$ is covered by an appropriate collection of sets $(\Delta_{\alpha})_{\alpha} \subset \mathbb{R}^2$ and one defines symplectic manifolds over them $(M_{\alpha}, F_{\alpha}:M_{\alpha} \to \Delta_{\alpha})$ using action\--angle coordinates. The second step is devoted to glueing these
models, using a symplectic glueing method which we describe below, in order to obtain a continuous
map 
$$
F \colon M-(\textup{neighborhoods of the focus\--focus fibers}) \longrightarrow \mathbb{R}^2,
$$
where by a focus\--focus fiber we mean a singular fiber containing one focus\--focus singularity.   
Next we attach the neighborhoods of the focus\--focus fibers, determined by the Taylor series. 
This glueing is delicate and requires an $\epsilon-\delta$ analysis near focus\--focus  fibers. There is remaining freedom when glueing in these fibers
which is controlled by the twisting index invariant.  
\item
 In this way, one obtains a continuous map $F \colon M \to \mathbb{R}^2$ 
which is ``essentially" the semitoric integrable system we want to construct, but has the problem that it
is singular. By singular we mean that it is non smooth (so we cannot yet speak of an integrable system), near the overlaps of the neighborhoods of the focus\--focus fibers 
(described by Eliasson's normal form) with the neighborhoods of the regular fibers (described by action\--angle
coordinates). 
\item
To overcome this smoothness problem of $F$ and produce a semitoric 
integrable system, one needs to suitably modify $F$ to make it smooth everywhere. 
This is a subtle analytic problem which takes the last step, and a significant part of the paper~\cite{PeVN11}.
\end{itemize}

To conclude, we discuss  how to do the symplectic glueing mentioned in the second
step above. The point is that this glueing works for continuous maps on symplectic manifolds,
and hence why above there needed to be a final step, to smoothen the continuous map which
the next result gives. The result can be formulated
 quite generally as follows (see Section~\cite[Section 3.4]{PeVN11} for a detailed formulation and proof).  
 Let us suppose that 
$(M_\alpha)_{\alpha\in A}$  are symplectic manifolds and 
$F_\alpha:M_\alpha\to V_\alpha\subset\mathbb{R}^n$ are proper continuous maps
where  $V_\alpha$  is open. For each
 $\alpha,\beta\in A$ assume that 
    $\varphi_{\alpha\beta} : F_{\alpha}^{-1}(D_{\alpha\beta})\to
    F_{\beta}^{-1}(D_{\alpha\beta})$ 
    is a symplectomorphism such that
        $ \varphi_{\alpha\beta}^* F_{\beta} = F_\alpha$ and
whenever
       $V_\alpha\cap V_{\beta}\cap V_\gamma\neq \emptyset$, one has that
    $\varphi_{\beta\gamma} \circ \varphi_{\alpha\beta}=
   \varphi_{\alpha\gamma}$.  Then one can prove that the smooth manifold
     $M$ obtained by glueing the symplectic manifolds
  $(M_\alpha)_{\alpha\in A}$
  by means of the transformations
  $(\varphi_{\alpha\beta})$ 
  is Hausdorff, paracompact, and symplectic, and that
  there exists a proper continuous map $F:M\to \bigcup_{\alpha\in A}
  V_\alpha\subset \mathbb{R}^n$
  such that $F_\alpha = F\circ y_\alpha$,
where  $y_\alpha:M_\alpha\hookrightarrow M$, $\alpha\in A$, 
are the natural inclusion symplectic embeddings.

\subsection{Semitoric systems with multiply pinched focus\--focus fibers in dimension $4$} \label{nss}

The classification of simple semitoric systems  in Sections~\ref{uni} and \ref{existence}
has been extended in~\cite{PPT19} to non\--simple systems. Hence 
the same singular fiber of $J$, and hence of $F$, may contain any finite
number $\lambda \in \mathbb{Z}^+$ of focus\--focus singularities; topologically such a fiber is
a torus with $\lambda$ pinched points, that is, a wedge of $\lambda$ spheres as shown in Figure~\ref{figure1}.

In this situation, a  
classification which extends the one explained earlier can be given in terms of 
the natural extensions of the invariants we have described.   
The complete symplectic invariant is a polygon (which may not be convex due to
a slight difference in the way it is constructed) with
$n$ marked points $p_1,\ldots, p_n$ inside of it, and for each such point
 $p_{\ell}$ there is assigned a label
$$p_{\ell} \rightsquigarrow \Big(k,  
\overbrace{\underbrace{\sum_{i,j} a^1_{ij}x^iy^j, \ldots, \sum_{i,j} a^{\lambda}_{ij}x^iy^j}_{\textup{One Taylor series per pinch in} \,\, F^{-1}(p_{\ell})}}^{=(\sum_{i,j} a^s_{ij}x^iy^j)_{s=1}^{\lambda}}\Big) \in \mathbb{Z} \times (\mathbb{R}[x,y])^{\lambda},$$ where there are as many Taylor
series corresponding to a particular point $p_{\ell}$ as there are pinched points  in the fiber 
$F^{-1}(p_{\ell})$. In a diagram:
$$
 \textup{semitoric system}  \underbrace{\rightsquigarrow}_{\textup{lost?}} 
  {\overbrace{\Delta,p_1,k_1,(\sum_{i,j} a^{1s}_{ij}x^iy^j)_{s=1}^{\lambda_1}
  \ldots
   p_n, k_n, (\sum_{i,j} a^{ns}_{ij}x^iy^j)_{s=1}^{\lambda_n}}^{\textup{symplectic invariants}}}  \underbrace{\rightsquigarrow}_{\textup{\S\ref{nss}}}  \textup{semitoric system}
$$

The same way as for simple semitoric systems, the answer to the question  in the first arrow of the diagram 
(whether data is lost in the process of finding the invariants)  is \emph{no} (up to isomorphisms).

These Taylor series were constructed in~\cite{PeTa}
as a generalization of the case when $\lambda=1$ from~\cite{VN03}. In the same way
as in Section~\ref{uni}, the integer $k$ and the Taylor series associated to each
point will vary with the point.

The precise classification statement~\cite[Theorem 4.10]{PPT19} is more involved
because  the twisting index invariant and the Taylor series invariants
 at the focus\--focus singularities are deeply connected if there exists at least one $\ell_0$ for
 which the number of focus\--focus singularities in $F^{-1}(p_{\ell_0})$ is $2$ or higher,
  and should not be considered
 as separate invariants. Instead a significant part of~\cite{PPT19}
 is devoted to describing a single invariant which incorporates the information of
 all the Taylor series and all the twisting indices simultaneously. 
 
Concrete examples of non simple semitoric systems appear in~\cite{DMH, HP}.
   
\subsection{Moduli spaces of simple semitoric  systems}

\subsubsection{Minimal models of toric and simple semitoric systems}
A toric  (respectively semitoric) integrable system $F \colon M \to \mathbb{R}^2$ is
\emph{minimal} if there is no toric (respectively semitoric) integrable system
 $F' \colon M' \to \mathbb{R}^2$ such that $F$ can be obtained from $F'$ by a blow up 
 respecting the toric (respectively semitoric) structure.

   There are three minimal models for toric integrable systems: their 
 associated fan corresponds to a square, a triangle, or a Hirzebruch trapezoid 
 (\cite[Theorem 8.2]{Oda} and \cite{Fulton}). 
 The symplectic geometry of minimal models of simple semitoric integrable
 systems was studied in~\cite{KPP-jgp}.  In this case there is an explicit
 list of minimal models in terms of a generalization of the fan; this generalization is called the \emph{helix}, and can be considered  a symplectic analogue of the fan of a nonsingular complete toric variety in algebraic geometry, that takes into account the effects of the monodromy near focus\--focus singularities.  By~\cite[Theorem 1.3]{KPP-jgp} 
 there are seven minimal models for simple semitoric integrable systems, corresponding to seven
 inequivalent helices.

 \subsubsection{A metric and a topology on the moduli spaces of toric and simple semitoric systems}
 One of the motivations to introduce the helix, and the other
 invariants of integrable systems (Taylor series, twisting index, semitoric polygon, etc),
 as well as the minimal models, is to shed light on the structure of the moduli
 spaces of the corresponding integrable systems which have these invariants. 

 The simplest cases to consider are probably the case of toric and simple
 semitoric integrable systems. 
  A first step in this direction was given in~\cite{PPRS}, where
 the case of toric  integrable systems on compact manifolds 
 is studied, and  a topology is defined on 
 the corresponding moduli space $\mathcal{M}_{\rm t}$. 
 
 The topology is constructed indirectly 
 by defining a  distance on the moduli space, and then
 considering the topology induced by this distance. In this paper the authors
 study a variety of topological properties of this space, in particular,
path-connectedness, compactness, and completeness. The construction of the distance
is related to the Duistermaat\--Heckman measure~\cite{DH}.   
This construction was generalized in~\cite{Palmer-moduli} to simple semitoric integrable systems, where the
author defines a metric on the moduli space of simple semitoric systems $\mathcal{M}_{\rm s}$, and studies
several of its properties.  The connectivity properties of this moduli space
were studied in~\cite[Section 6.3]{KPP-sigma}.

\subsubsection{Functions on moduli spaces of toric and simple semitoric systems}
One of the useful consequences of having topologies on these moduli spaces 
is that now one is able to quantify the variation of certain natural functions
$\mathcal{M}_{\rm t} \to \mathbb{R}$ and $\mathcal{M}_{\rm s} \to \mathbb{R}$ defined on them. 
For instance, on $\mathcal{M}_{\rm t}$ with the topology  mentioned earlier, one can define the
so called \emph{toric packing capacity}  $\mathcal{M}_{\rm t} \to \mathbb{R}$
as follows.

 First, let ${\rm B}^{2n}(r)\subset \mathbb{R}^{2n} \simeq \mathbb{C}^n$ be the ball
of radius $r$ centered at the origin. It comes endowed with
a natural rotational action of the $n$\--torus $\mathbb{T}^n$ component by component of $\mathbb{C}^n$. 
With this notation in mind,  a \emph{ball packing} $P$ of $(M,\omega)$ is any disjoint
union of symplectically embedded balls ${\rm B}^{2n}(r) \hookrightarrow M$ into
$M$, of any possibly varying radii $r>0$. A \emph{toric ball packing}~\cite{Pe06,Pe07}
of  a toric system $F \colon M \to \mathbb{R}^n$ is a ball packing of $(M,\omega)$
such that  the symplectic embeddings of balls are \emph{equivariant} with respect to the natural
$n$\--torus action on ${\rm B}^{2n}(r)$ and the $n$\--torus action 
induced on $M$ by $F$. The \emph{volume of}
$P$ is defined to be the sum of the volumes of the balls with respect
to the volume form $\omega^n$; we denote it by ${\rm vol}(P)$. With these
definitions in place we may define a function $ c \colon \mathcal{M}_{\rm t} \to \mathbb{R}$ by
assigning to each toric integrable system the number
\begin{eqnarray}
c(F)= \Big(  \frac{1}{{\rm vol} ({\rm B}^{2n})}
 \sup \{{\rm vol}(P)\,|\, P \,\,\textup{is a toric ball packing of}\,\, F\}   \Big)^{\frac{1}{2n}}. \nonumber
\end{eqnarray}

An analogous function can be defined on $\mathcal{M}_{\rm s}$.  These functions are examples
of what are called in~\cite{FPP} \emph{$G$\--equivariant symplectic
capacities}, which generalize to the equivariant setting (here $G$
is any Lie group) the usual notion of symplectic capacity.

The continuity of these (and closely related) functions was studied in~\cite{FP}
and then in~\cite{FPP}. Indeed, by \cite[Theorem A]{FP} and \cite[Theorem 1.2]{FPP} the function
 $ c \colon \mathcal{M}_{\rm t} \to \mathbb{R}$ is everywhere discontinuous
and the restriction to the subspace of toric integrable systems
with exactly $N$ points fixed by the induced $n$\--torus action is continuous for any
choice of $N>0$. Also in  \cite[Theorem 1.2]{FPP} an analogous continuity statement
is given for semitoric integrable systems.

\subsection{Other classifications or related results}

The symplectic theory of toric integrable systems discussed in Section~\ref{toricsection}, and the classifications of the Fomenko school~\cite{BS}, were two of the motivations for pursuing the global symplectic 
classification of simple semitoric systems and its generalizations which we have discussed in
Sections~\ref{uni}, \ref{existence}, and~\ref{nss}. 

The classifications we have presented 
have relations to almost toric systems and systems with semitoric 
features as in~\cite{HSSS2, ls, s, Wa1, Wa2}, and also 
the theory of Hamiltonian $S^1$\--spaces~\cite{HSSS, K}.  

There is also recent work on the so called toric\--focus systems~\cite{RWZ} and systems which have
semitoric features but also include  singularities with hyperbolic components~\cite{DP}. In
the article~\cite{Tang19} there is an application of the ideas in the proof of this classification to study
symplectic forms on noncompact manifolds.

\section{Inverse spectral geometry of quantum toric or semitoric systems} \label{spec1}

A quantum integrable system  is given by a collection of semiclassical
commuting self\--adjoint operators $P_1,\ldots, P_n$ whose principal symbols form a classical 
integrable system. 
Some well known examples of quantum integrable systems are the quantum spherical pendulum, discussed by
Cushman and Duistermaat in \cite{[ab10]} and the ``Champagne bottle"
\cite{[ab6]}. Quantum integrable systems given by Berezin\--Toeplitz quantization are common
in the physics literature. An  example is the coupled angular momenta \cite{LFPe19, SZ}; see 
also \cite[Section 8.3]{PPVN} for a proof that it is a Berezin\--Toeplitz
system.  

The study of quantum integrable systems such as these examples fits into a general framework
of ideas, which we briefly discuss next, which uses a combination
of techniques from symplectic geometry and microlocal analysis to go back and forth
between classical and quantum systems.

\subsection{The  inverse spectral problem for quantum integrable systems}

The inverse problem we are going to discuss next belongs to a class of semiclassical inverse spectral
questions which has been the focus of intense attention in recent years~\cite{CoGu2011, Gu95, Ha13,  LF14, 
 Pe-bbms, PPVN, VN2011}.  The problem goes back to pioneer works of B\'erard~\cite{Be76},
Br\"uning\--Heintze~\cite{BH84}, Colin de Verdi{\`e}re~\cite{CdV,
  CdV2}, Duistermaat\--Guillemin~\cite{DuGu1975}, and
Guillemin\--Sternberg \cite{GuSt}, in the 1970s/1980s, 
and is closely related to inverse problems that are not directly
semiclassical but use similar microlocal methods~\cite{Ze98, Ze04}.

The following statement  which corresponds to~\cite[Conjecture 9.1]{BAMS11},
concerns only simple semitoric systems; we keep the original formulation, including 
the notation used therein: ``a semitoric system $J,H$ is determined up to symplectic
   equivalence by its semiclassical joint spectrum as $\hbar \to 0$.
   From any such spectrum one can construct explicitly the associated
   semitoric system, i.e.  the set of points in $\mathbb{R}^2$ where
   on the $x$\--basis we have the eigenvalues of $\hat{J}$, and on the
   vertical axis the eigenvalues of $\hat{H}$ restricted to the
   $\lambda$\--eigenspace of $\hat{J}$".   
     A motivation to develop the symplectic geometric results of the previous sections
has been to shed light on this conjecture. 

The proof strategy, outlined in~\cite{PeVN12-dcds}, 
 consists of detecting the symplectic invariants (Section~\ref{uni}) in the joint spectrum, and use them
 to construct the associated classical system. The semiclassical inverse spectral problem can be
posed more generally as the question of how much information about the associated classical
system can be obtained from  the knowledge of the semiclassical spectrum
of a quantum integrable system:
$$
\overbrace{(P^1_{\hbar})_{\hbar \in I}, \ldots, (P^n_{\hbar})_{\hbar \in I}}^{\textup{quantum integrable system}}
\rightsquigarrow
 \overbrace{(X_{\hbar})_{\hbar \in I} \subset \mathbb{R}^n}^{\textup{one joint spectrum for each}\,\, \hbar}
 \underbrace{\rightsquigarrow}_{\textup{possible?}} 
 \overbrace{f_1,\ldots, f_n}^{\textup{principal symbols of}\,\,(P^1_{\hbar})_{\hbar \in I}, \ldots, (P^n_{\hbar})_{\hbar \in I} },
 $$ 
 where $f_1,\ldots, f_n$ stand for the principal symbols of the operators on the left hand side.

\subsection{Progress towards the inverse spectral problem in the simple semitoric case}

Next we will briefly discuss how the general approach we mentioned earlier can be implemented
to shed light on the inverse spectral conjecture for simple semitoric systems. We start by presenting
the basic language of spectral geometry.

\subsubsection{The joint spectrum}
\label{sec:semicl-oper}

Let $(M,\omega)$ be a connected symplectic manifold of dimension $4$ 
and let $I$ be any subset of $(0,1]$ which has  as an accumulation point at $0$. 
For any complex Hilbert space $\mathcal{H}$ we denote by
$\mathcal{L}(\mathcal{H})$ the set of linear, possibly unbounded,
self-adjoint operators on $\mathcal{H}$ with a dense domain.

Consider a sequence of Hilbert spaces $(\mathcal{H}_{\hbar})_{\hbar \in I}$.
A space $\Psi$ of \emph{semiclassical operators} is a subspace of
$\prod_{\hbar \in I} \mathcal{L}(\mathcal{H}_{\hbar})$, which contains the
identity and is equipped with a weakly positive principal symbol map,
which is a (normalized) $\mathbb{R}$-linear map $\sigma \colon \Psi \to
\mathcal{C}^{\infty}(M,\mathbb{R})$ satisfying a \emph{product formula} (if $P, Q$ are in $\Psi$ and if $P\circ Q$ is well defined and is in  the space $\Psi$ then  we have an equality $\sigma(P\circ Q) = \sigma(P)\sigma(Q)$) and a \emph{weak positivity condition} (if $\sigma(P)\geqslant 0$  there is a  function $\hbar\mapsto\epsilon(\hbar)$ tending to zero as $\hbar\to 0$ and such
  that one has $P\geqslant -\epsilon(\hbar)$ for all $\hbar\in I$).  
  
Two famous examples of such semiclassical operators are
semiclassical pseudodifferential operators~\cite{DiSj99, Zw2012} and
semiclassical (or Berezin\--)Toeplitz operators~\cite{BG81, BorPauUri,Ch2003b, Ch2006b,
  Ch2006a, Ch2007,MaMar2008,Schli2010}.

  If $P=(P_{\hbar})_{\hbar \in I} \in \Psi$, the image 
  $\sigma(P)$ is
called the \emph{principal symbol of $P$}.   The principal symbol plays a fundamental
role in the formulation of the inverse spectral problem for quantum integrable
systems which we will give shortly.

Let $P=(P_{\hbar})_{\hbar \in I}$ and $Q=(Q_{\hbar})_{\hbar \in I}$ be
semiclassical operators on $(\mathcal{H}_\hbar)_{\hbar \in I}$. We say that $P$ and $Q$ \emph{commute} if
for each $\hbar \in I$ the operators $P_{\hbar}$ and $Q_{\hbar}$ commute 
(note that in the case of unbounded self\--adjoint operators, by definition 
this means that their projector\--valued spectral measures commute). In this case one
can define for each fixed value of the parameter $\hbar$ the so called 
 \emph{joint spectrum} of $(P_{\hbar},Q_{\hbar})$ as the support of
 the joint spectral measure. We denote it by
 ${\rm JointSpec}(P_{\hbar},\,Q_{\hbar})$. 
 
 If the Hilbert space $\mathcal{H}_\hbar$ is
 finite dimensional, or more generally, if the joint spectrum is
 discrete, 
 $$
 {\rm JointSpec}(P_{\hbar},Q_{\hbar})=\Big\{(\lambda_1,\lambda_2)\in
 \mathbb{R}^2\,\, |\,\, \exists v\neq 0,\,\, P_{\hbar} v = \lambda_1
 v,\,\,Q_{\hbar} v = \lambda_2 v \Big\}.
 $$
 The \emph{joint spectrum} of $(P,Q)$ is, by definition, the collection of all 
 ${\rm JointSpec}(P_{\hbar},Q_{\hbar}), \hbar \in I$, and  it is denoted by
 ${\rm JointSpec}(P,\,Q).$ For convenience we also
 view  ${\rm JointSpec}(P,\,Q)$ as a set depending on $\hbar$. 
 
 All of these
 definitions extend easily to $2n$\--dimensional manifolds and 
 collections of commuting semiclassical operators $P_1,\ldots, P_n$.

\subsubsection{Bohr\--Sommerfeld rules}

The abstract inverse spectral result for simple semitoric systems we are going to present in the next section
requires that the operators follow some well known properties, known as the \emph{Bohr\--Sommerfeld
rules}. These are easy to state, and we are going to formulate them next.  They are  
 also known to hold for integrable systems of  pseudodifferential operators~\cite{charbonnel,vungoc-focus},  where  the manifold $M$ is a cotangent bundle,  and for integrable systems of Berezin\--Toeplitz
operators on prequantizable compact symplectic
manifolds~\cite{charles-quasimodes}.

 Let $F=(f_1,f_2) \colon (M,\omega) \to
  \mathbb{R}^2$  be an integrable system on a connected symplectic manifold of dimension $4$. Let $P$ and $Q$ be commuting
  semiclassical operators with principal symbols $f_1,f_2 \colon M \to
  \mathbb{R}$. We say
  that ${\rm JointSpec}(P,\,Q)$ \emph{satisfies the Bohr\--Sommerfeld
    rules} 
  if for every regular value $c$ of $F$ we can find a ball
  ${\rm B}(c,\epsilon_c)$ centered at $c$ (of some radius as small as necessary), such that,    on
    ${\rm
      B}(c,\epsilon_c)$, we have
  $$
 {\rm JointSpec}(P,Q) = g_{\hbar}(2\pi \hbar
    \mathbb{Z}^2\cap D) + \mathcal{O}(\hbar^2)
    $$
  with $g_{\hbar}=g_0+\hbar g_1$, where $g_0,g_1$ are smooth maps
  defined on a bounded open set $D\subset\mathbb{R}^2$, $g_0$ is a
  diffeomorphism into its image, $c \in g_0(D)$ and the components of
  $g_0^{-1}=(\mathcal{A}_1,\,\mathcal{A}_2)$ are such that $(\mathcal{A}_1\circ F, \mathcal{A}_2\circ F)$ form a basis of action
  variables in the sense that the two Hamiltonians in 
  action\--angle coordinates correspond to $(\xi_1,\xi_2)$ in Section~\ref{acc}. 
  
  The displayed
  equation is a precise statement if one uses the Hausdorff distance ${\rm d}_H$ to compare how
  far a set is from another; with this in mind, if for instance  $(A_{\hbar})_{\hbar \in I}$ and
$(B_{\hbar})_{\hbar \in I}$ are sequences of uniformly bounded subsets
of $\mathbb{R}^2$, then the notation $A_{\hbar} = B_{\hbar} +
\mathcal{O}(\hbar^{N}) $ means that there is a constant $C>0$ such that
$
{\rm d}_H(A_{\hbar},\,B_{\hbar})\leqslant C\hbar^{N}$ for all $\hbar \in
I$, see~\cite[Sections 2.6 and 2.7]{LFPeVN16}.

\subsubsection{Results for all invariants with the exception of the twisting indices} \label{ss1}

A \emph{quantum semitoric integrable system} $(P,Q)$ is given by two semiclassical
commuting self\--adjoint operators whose principal symbols form a classical semitoric
integrable system. The quantum semitoric integrable system $(P,Q)$ is \emph{simple}
if the corresponding classical semitoric integrable system is also simple.

The main result known for simple quantum semitoric systems~\cite[Theorem A]{LFPeVN16} can be stated as follows. Let $(P,Q)$ be a simple quantum  semitoric system on $M$ for which the
Bohr\--Sommerfeld rules hold.  Then from the knowledge of ${\rm JointSpec}(P, Q)+\mathcal{O}(\hbar^2),$ one can recover all of the invariants of the symplectic classification in Sections~\ref{uni} and \ref{existence}, 
except for the twisting index $k$, namely:
\begin{itemize}
\item
 the convex polygon, 
 \item
the marked points $p_1,\ldots,p_n$ in the interior of the polygon, and
\item
 for each point above the Taylor series $\sum_{i,j} a_{ij}\,x^iy^j$.
 \end{itemize}

One needs to combine microlocal analytic techniques and symplectic geometry in order to
recover from the joint spectrum the polygon, the interior marked points
$p_1,\ldots,p_n$ and for each such point the label $\sum_{i,j} a_{ij}\,x^iy^j$. For this
the crucial part is that we have a semiclassical spectrum, so we know the variation
with respect to $\hbar$ for a sequence of values converging to $0$. This is analyzed
with the help of the Bohr\--Sommerfeld rules.

The image  of the joint map of principal symbols 
$(f_1,f_2)(M)$ is the limit as $\hbar \to 0$ of the joint spectrum.
To recover the polygon invariant and the marked points is a more subtle task because
here it matters the position of the images of the focus\--focus singularities inside of $F(M)$,
reflected in how and where the joint eigenvalues (the elements of the spectrum) 
concentrate, as $\hbar \to 0$. This step consists essentially of recovering the singular
affine structure induced by the fibration given by the integrable
system near the focus\--focus singularities.  For this  we need
to use Eliasson's linearization theorem and
the work on singular affine structures induced by 
semitoric systems~\cite{PeVN09, PeVN11, VN07},
among other ingredients. In this way one recovers the
integral affine structure from the spectrum modulo $\mathcal{O}(\hbar^2)$. 

In order to recover $\sum_{i,j}a_{i,j}x^i y^j$ at the focus\--focus singularities
one needs to use Bohr\--Sommerfeld rules, semiclassical Fourier transforms and 
the lemma of non stationary phase. To recover the explicit form of $\sum_{i,j}a_{i,j}x^i y^j$ 
is essential to understand its construction in terms of the
period maps $\tau_1$ and $\tau_2$ introduced in Section~\ref{uni}.

\subsection{Solution to the inverse spectral problem in the compact toric case} \label{ts}

The inverse spectral problem problem was
solved in the positive in~\cite{CPVN13}, for quantum toric
integrable systems given by Berezing\--Toeplitz operators on compact manifolds 
(prior to the developments in the semitoric case in Section~\ref{ss1}). Recall
from Section~\ref{toricsection} that we always assume that toric integrable systems
take place on compact connected manifolds, but the dimension of
the manifold can be arbitrary.

In analogy with the semitoric case which we have just 
discussed, by a \emph{quantum toric  integrable system} we mean a quantum integrable system 
$P_1, \ldots, P_n$ on a symplectic $2n$\--dimensional manifold $(M,\, \omega)$   such that the
principal symbols $f_1,\ldots, f_n$ of   $P_1, \ldots, P_n$ form a toric integrable system.

The proof strategy in this case is the same as in the semitoric case,
with the simplification that there is only one symplectic invariant,
the image $F(M) \subset \mathbb{R}^n$, as shown in Section~\ref{delzant} (see Figure~\ref{figure2}).

The paper~\cite{CPVN13} contains a detailed analysis of the semiclassical
spectral theory of toric integrable systems, and as a consequence of it, it is shown the
semiclassical joint spectrum of $P_1, \ldots, P_n$ converges to the Delzant
polytope $\Delta=F(M)$. This is why in this case the dimension of $M$ is not
essential (as it was in the semitoric case, where our results have been restricted
to dimension $4$); what is essential is that $M$ is compact (there is a recent extension
of Delzant's theorem to noncompact manifolds~\cite{KL2015} which may
be useful in generalizing~\cite{CPVN13} to the noncompact case). This result was motivated
and preceded by~\cite{CdV, CdV2}.

As a consequence of the semiclassical spectral theory developed therein 
it was shown that all toric integrable systems on compact connected manifolds
can be quantized~\cite[Theorem 1.4]{CPVN13}. That is, if $F:=(f_1,\ldots, f_n) \colon M \to \mathbb{R}^n$ is
a toric integrable system on a compact manifold then there exists a quantum
toric integrable system $\psi_1,\ldots, \psi_n$  with associated principal symbols
$f_1,\, \ldots,\, f_n$.

\subsection{Other results}

 The result in Section~\ref{ss1}  gives no information
on the twisting index $k$ for any of the points. Nonetheless it sheds light on  the inverse
spectral conjecture. 

An inverse spectral result was given in~\cite{PeVN14-cmp} for a neighborhood of
a singular fiber of an integrable system containing exactly one focus\--focus singularity $m$.  
It was shown 
therein that the joint spectrum of the quantum integrable system near the singular fiber 
determines the Taylor series invariant at $m$, and hence
determines symplectically a neighborhood of the fiber up to isomorphisms of integrable systems. 
This result was needed for the global result described in Section~\ref{ss1}.

\begin{figure}[h]
\centering 
\includegraphics[width=0.9\textwidth]{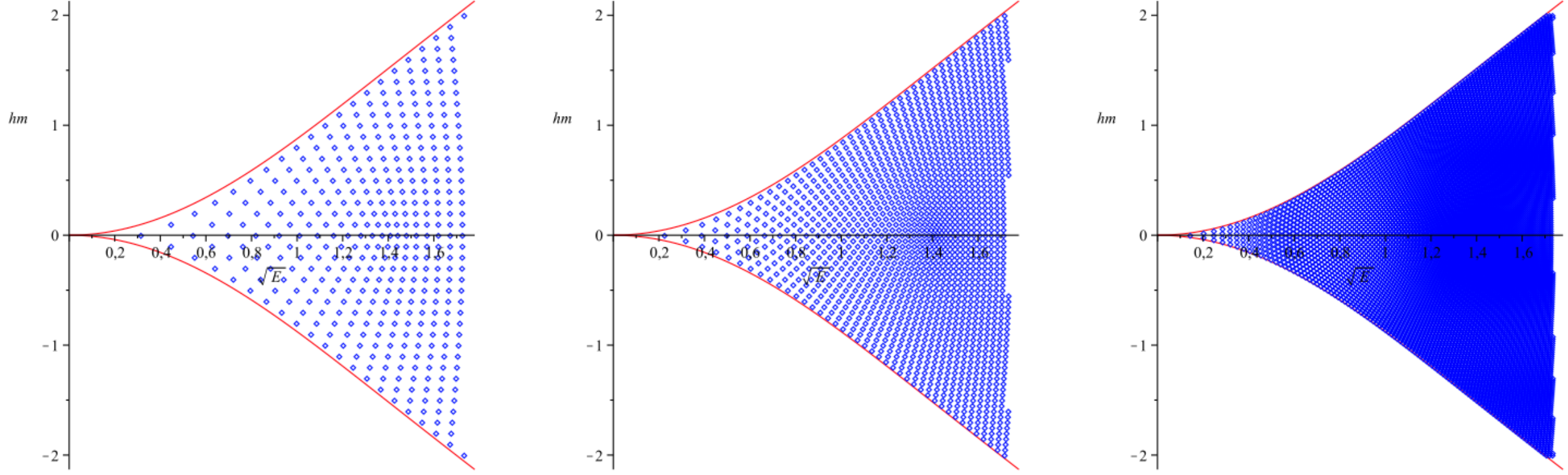}
 \caption{Joint spectrum of normalized quantum spherical pendulum as $\hbar \to 0$.}
\label{figure4}
\end{figure}

For general classes of self adjoint pseudodifferential operators and for Berezin\--Toeplitz
operators on compact manifolds it was shown in~\cite{PeVN-montesinos} that
$$
\lim_{\hbar \rightarrow 0} \textup{JointSpec}(P_1,\dots, P_d)  = F(M) \subset \mathbb{R}^d,
$$ 
where $F$ is the joint map formed by the principal symbols $P_1,\ldots, P_d$. This convergence is illustrated 
in Figure~\ref{figure4} for the 
quantum spherical pendulum which shows, as originally illustrated
in~\cite[Figure 1]{PeVN-montesinos}, that 
as the parameter $\hbar$ approaches the accumulation point $0$ the semiclassical joint spectrum  fills
 the inside of the curve in red color. This curve is the boundary of the image of the
 joint map of principal symbols $F=(f_1,f_2) \colon M \to \mathbb{R}^2$.

There are several related results  in~\cite{PPVN} which 
illustrate the general phenomena of convergence of the semiclassical spectrum, in the semiclassical 
limit, to the joint image of the principal symbols in $\mathbb{R}^n$ (which is often called
the \emph{classical spectrum}). These results were extended to the case of unitary operators 
in~\cite{LFPe19-jmap}.

   The articles~\cite{BMMT, PeVN12-dcds} contain a number of open questions and 
   problems about classical and quantum integrable systems as well as some strategies
   to approach them.

\subsection{A question}

If the  inverse spectral conjecture~\cite[Conjecture 9.1]{BAMS11} holds as stated for simple semitoric
systems,  it would be interesting to see if a similar statement can be given 
for non\--simple semitoric systems (in Section~\ref{nss}), and we can ask the question (essentially from~\cite[Section 1]{PPT19}): 

\medskip

\paragraph{{\bf Question}} \emph{Are there two  non\--simple  quantum semitoric integrable systems $(P,Q)$ and $(P',Q')$
with the same joint spectrum {\rm(}modulo $\mathcal{O}(\hbar^2)${\rm)} and such that
the associated classical systems of principal symbols $(f_1,f_2)$ and $(f'_1,f'_2)$ are not isomorphic?}

\medskip

It would be interesting to study this question for  the integrable systems in~\cite{DMH}. 

A solution to this question would be helpful in understanding to what extent 
the structure or composition of the singular set of a classical integrable system $(f_1,f_2)$ can be read 
off from the spectrum the quantum integrable system $(P,Q)$ with principal symbols given by
$(f_1,f_2)$. We saw in Section~\ref{ts} that all of the symplectic geometry of
 toric integrable systems on compact connected symplectic manifolds 
 can be read off from the spectrum, up to isomorphisms.
      
\subsubsection*{Acknowledgements}
 
 \medskip
 
I am very grateful to the Departamento de \'Algebra, Geometr\'ia y Topolog\'ia at
the Universidad Complutense de Madrid for the excellent hospitality during a visit in
 December of 2019 and January 2020, when part of this paper was written.  I would also like
to thank Yohann Le Floch, Joseph Palmer, and Xiudi Tang, 
for  comments on a preliminary version of this paper, and Xiudi Tang for help drawing the figures.

This paper is dedicated to Professor J.J. Duistermaat.  
I was fortunate to collaborate with him on four
papers on symplectic group actions. His work on symplectic geometry and analysis continues to 
influence my own research
 on symplectic geometry of integrable systems and group actions.   The articles~\cite{NAMS, AP11} discuss some of Professor Duistermaat's contributions.

\noindent
\\
{\'A}lvaro Pelayo\\
Department of Mathematics\\
University of California, San Diego\\
9500 Gilman Drive  \# 0112\\
La Jolla, CA 92093-0112, USA\\
{\em E\--mail}: \texttt{alpelayo@math.ucsd.edu}

\end{document}